\documentclass[11pt]{article}

\usepackage{amsmath,amssymb}
\usepackage[mathscr]{eucal}
\usepackage{theorem}

\newtheorem{theorem}{Theorem}[section]
\newtheorem{proposition}[theorem]{Proposition}
\newtheorem{lemma}[theorem]{Lemma}
\newtheorem{corollary}[theorem]{Corollary}
{\theorembodyfont{\rmfamily}
\newtheorem{remark}[theorem]{\rm R e m a r k \ }}

\makeatletter
\@addtoreset{equation}{section}
\makeatother

\begin{document}

\large

\title{\Huge\bf Function Theory on a q-Analog of Complex Hyperbolic Space}
\author{\Large O. Bershtein \and \Large S. Sinel'shchikov}
\date{Mathematics Division, B. Verkin Institute for Low Temperature Physics
and Engineering,
\\ National Academy of Sciences of Ukraine
\\ 47 Lenin Ave., Kharkov 61103, Ukraine}

\maketitle

\begin{abstract}
    This work deals with function theory on quantum complex hyperbolic
spaces. The principal notions are expounded. We obtain explicit formulas for invariant integrals on `finite'
functions on a quantum hyperbolic space and on the associated quantum isotropic cone. Also we establish
principal series of $U_q \mathfrak{su}_{n,m}$-modules related to this cone, and obtain the necessary
conditions for those modules to be equivalent.
\end{abstract}

\section{Introduction}

Let us consider the group $SU_{n,m}$ of pseudo-unitary $(n+m) \times
(n+m)$-matrices that preserve the following form in $\mathbb C^{n+m}$:
\begin{equation*}
[x,y]=-x_1 \bar{y}_1-\ldots -x_n \bar{y}_n + x_{n+1}\bar{y}_{n+1}+\ldots
+x_{n+m}\bar{y}_{n+m}.
\end{equation*}
Then one can also consider the manifold $\widehat{\mathscr{H}}_{n,m}=\{x
\in \mathbb C^{n+m}|[x,x]>0\}$ and its projectivization
$\mathscr{H}_{n,m}$. The latter manifold is isomorphic to the homogeneous
space $SU_{n,m}/S(U_{n,m-1}\times U_1)$, a complex hyperbolic space. There
is a vast literature devoted to the study of these pseudo-Hermitian spaces
of rank 1, in particular harmonic analysis on those (see J.Faraut
\cite{Faraut}, V.Molchanov \cite{Molch1,Molch2}, G.van Dijk and Yu.Sharshov
\cite{Dijk-Sh}).

In this paper we establish basic notions in the theory of quantum
pseudo-Hermitian spaces. These objects initially appear in the work of
Reshetikhin, Faddeev and Takhtadjan \cite{RTF}. Later on the development of
the theory of quantum bounded symmetric domains and quantum analogs of
representation theory of noncompact real Lie groups made it clear that the
above objects really worth studying. For example, the Penrose transform of
the quantum matrix ball of rank 2 leads to a quantum analog of the complex
hyperbolic space in $\mathbb C^4$, see \cite{Penrose}.

We introduce a background of the function theory on quantum analogs of complex hyperbolic spaces
$\mathscr{H}_{n,m}$ and of the related isotropic cones $\Xi_{n,m}=\{x \in \mathbb C^{n+m}|[x,x]=0\}$. We
establish some special `spaces of functions with compact support' (called finite functions, for short) and
endow these noncommutative algebras with faithful representations. Then we introduce integrals on the spaces
of finite functions and prove their invariance under the action of quantum universal enveloping algebra $U_q
\mathfrak{su}_{n,m}$. Finally, we introduce a quantum analog of the principal (unitary) series of $U_q
\mathfrak{su}_{n,m}$-modules related to a quantum analog of the cone $\Xi$. For these modules we establish
the necessary conditions for the equivalence.

This project started out as joint work with L.~Vaksman and D.~Shklyarov. We are grateful to
both of them for helpful discussions and drafts with preliminary definitions and
computations.

\medskip

\section{Preliminaries}

Let $q\in(0,1)$. The Hopf algebra $U_q\mathfrak{sl}_{N}$ is given by its
generators $K_i$, $K_i^{-1}$, $E_i$, $F_i$, $i=1,2,\ldots,N-1$, and the
relations:
$$K_iK_j=K_jK_i,\qquad K_iK_i^{-1}=K_i^{-1}K_i=1,$$
$$K_iE_i=q^{2}E_iK_i,\qquad K_iF_i=q^{-2}F_iK_i,$$
$$K_iE_j=q^{-1}E_jK_i,\qquad K_iF_j=qF_jK_i,\qquad |i-j|=1,$$
$$E_iF_j-F_jE_i=\delta_{ij}\frac{K_i-K_i^{-1}}{q-q^{-1}},$$
$$E_i^2E_j-(q+q^{-1})E_iE_jE_i+E_jE_i^2=0,\qquad|i-j|=1,$$
$$F_i^2F_j-(q+q^{-1})F_iF_jF_i+F_jF_i^2=0,\qquad|i-j|=1,$$
$$[E_i,E_j]=[F_i,F_j]=0,\qquad|i-j|\ne 1.$$
The comultiplication $\Delta$, the antipode $S$, and the counit
$\varepsilon$ are defined on the generators by
$$
\Delta(E_i)=E_i\otimes 1+K_i\otimes E_i,\quad
\Delta(F_i)=F_i\otimes K_i^{-1}+1\otimes F_i,\quad
\Delta(K_i)=K_i\otimes K_i,
$$
$$S(E_i)=-K_i^{-1}E_i,\qquad S(F_i)=-F_iK_i,\qquad S(K_i)=K_i^{-1},$$
$$\varepsilon(E_i)=\varepsilon(F_i)=0,\qquad\varepsilon(K_i)=1,$$
see \cite[Chapter 4]{Jant}.

We need also the Hopf algebra $\mathbb{C}[SL_{N}]_q$ of matrix elements of
finite dimensional weight $U_q\mathfrak{sl}_{N}$-modules. Recall that
$\mathbb{C}[SL_{N}]_q$ can be defined by the generators $t_{ij}$,
$i,j=1,...,N$, (the matrix elements of the vector representation in a
weight basis) and the relations
\begin{align*}
&t_{ij'}t_{ij''}=qt_{ij''}t_{ij'},\qquad &j'<j'',
\\ &t_{i'j}t_{i''j}=qt_{i''j}t_{i'j},\qquad &i'<i'',
\\ &t_{ij}t_{i'j'}=t_{i'j'}t_{ij},\qquad &i<i'\;\&\;j>j',
\\ &t_{ij}t_{i'j'}=t_{i'j'}t_{ij}+(q-q^{-1})t_{ij'}t_{i'j},\qquad
& i<i'\;\&\; j<j',
\end{align*}
together with one more relation
\begin{equation*}\label{det_1}
\det\nolimits_q\mathbf{t}=1,
\end{equation*}
where $\det\nolimits_q\mathbf{t}$ is a $q$-determinant of the matrix
$\mathbf{t}=(t_{ij})_{i,j=1,...,N}$:
\begin{equation*}
\det\nolimits_q\mathbf{t}=\sum\limits_{s\in S_{N}}(-q)^{l(s)}t_{1 s(1)}
t_{2s(2)}\ldots t_{Ns(N)},
\end{equation*}
with $l(s)=\mathrm{card}\{(i,j)|i<j\;\&\;s(i)>s(j)\}$.

Let also $U_q\mathfrak{su}_{n,m}$, $m+n=N$, denotes the Hopf $*$-algebra $(U_q \mathfrak{sl}_N,*)$ given by
$$
(K_j^{\pm 1})^*=K_j^{\pm 1},\qquad E_j^*=
\begin{cases}
K_jF_j,& j \ne n,
\\ -K_jF_j,& j=n,
\end{cases}\qquad F_j^*=
\begin{cases}
E_jK_j^{-1},& j \ne n,
\\ -E_jK_j^{-1},& j=n,
\end{cases}
$$
with $j=1,\ldots,N-1$ \cite{RTF, polmat}.

\medskip

\section{\boldmath $*$-Algebra
$\operatorname{Pol}\left(\mathscr{H}_{n,m}\right)_{q}$}

Let $m,n\in\mathbb{N}$, $m\ge 2$, and $N\overset{\mathrm{def}}{=}n+m$. Recall that the
classical complex hyperbolic space $\mathscr{H}_{n,m}$ can be obtained by
projectivization of the domain
$$
\widehat{\mathscr{H}}_{n,m}=\left\{(t_1,\ldots,t_N)\in\mathbb{C}^N\left|\:
-\sum^n_{j=1}|t_j|^2+\sum^N_{j=n+1}|t_{j}|^2>0\right.\right\}.
$$

Now we pass from the classical case $q=1$ to the quantum case $0<q<1$. Let us consider the well known
\cite{RTF} $q$-analog of the pseudo-Hermitian spaces. Let
$\operatorname{Pol}\left(\widehat{\mathscr{H}}_{n,m}\right)_q$ denotes the unital $*$-algebra with the
generators $t_1,t_{2},\ldots,t_N$ and the commutation relations as follows:
\begin{equation}\label{t_it_j}
\begin{aligned}
t_it_j &= qt_jt_i,\qquad i<j
\\ t_it_j^* &= qt_j^*t_i,\qquad i\ne j
\\ t_it_i^* &= t_i^*t_i+(q^{-2}-1)\sum_{k=i+1}^Nt_kt_k^*,\qquad i>n
\\ t_it_i^* &= t_i^*t_i+(q^{-2}-1)\sum_{k=i+1}^nt_kt_k^*-
(q^{-2}-1)\sum_{k=n+1}^Nt_kt_k^*,\qquad i\le n.
\end{aligned}
\end{equation}
It is important to note that
$$c=-\sum_{j=1}^nt_jt_j^*+\sum_{j=n+1}^Nt_jt_j^*$$
is central in $\operatorname{Pol}\left(\widehat{\mathscr{H}}_{n,m}\right)_q$. Moreover,
$c$ is not a zero divisor in
$\operatorname{Pol}\left(\widehat{\mathscr{H}}_{n,m}\right)_q$. This allows one to embed
the $*$-algebra $\operatorname{Pol}\left(\widehat{\mathscr{H}}_{n,m}\right)_q$ into its
localization $\operatorname{Pol}\left(\widehat{\mathscr{H}}_{n,m}\right)_{q,c}$ with
respect to the multiplicative system $c^{\mathbb{N}}$.

The $*$-algebra $\operatorname{Pol}\left(\widehat{\mathscr{H}}_{n,m}\right)_{q,c}$ admits
the following bigrading:
$$\deg t_{j}=(1,0),\qquad\deg t_j^*=(0,1),\qquad j=1,2\ldots,N.$$
Introduce the notation
$$
\operatorname{Pol}(\mathscr{H}_{n,m})_q=\left\{\left.f\in
\operatorname{Pol}\left(\widehat{\mathscr{H}}_{n,m}\right)_{q,c} \right|\:\deg
f=(0,0)\right\}.
$$
This $*$-algebra $\operatorname{Pol}(\mathscr{H}_{n,m})_q$ will be called the algebra of
regular functions on the quantum hyperbolic space.

We are going to endow the $*$-algebra $\operatorname{Pol}(\mathscr{H}_{n,m})_q$ with a
structure of $U_{q}\mathfrak{su}_{n,m}$-module algebra \cite{Ch-P}. For this purpose, we
embed it into the $U_{q}\mathfrak{su}_{n,m}$-module $*$-algebra
$\operatorname{Pol}\left(\widetilde{X}\right)_q$ of `regular functions on the quantum
principal homogeneous space' constructed in \cite{polmat}.

Recall that
$\operatorname{Pol}\left(\widetilde{X}\right)_q\overset{\mathrm{def}}{=}
(\mathbb{C}[SL_N]_q,*)$, with $\mathbb{C}[SL_N]_q$ being the well-known
algebra of regular functions on the quantum group $SL_N$, and the
involution $*$ being defined by
$$
t_{ij}^*=\mathrm{sign}[(i-m-1/2)(n-j+1/2)](-q)^{j-i}\det\nolimits_qT_{ij}.
$$
Here $\det_q$ is the quantum determinant \cite{Ch-P}, and the matrix $T_{ij}$ is derived from the matrix
$T=(t_{kl})$ by discarding its $i$'s row and $j$'s column.

It follows from $\det_qT=1$ that
\begin{equation*}\label{bhr}
-\sum_{j=1}^nt_{1j}t_{1j}^*+\sum_{j=n+1}^Nt_{1j}t_{1j}^*=1.
\end{equation*}
Thus the map $J:t_j\mapsto t_{1j}$, $j=1,2,\ldots,N$, admits a unique extension to a
homomorphism of $*$-algebras
$J:\operatorname{Pol}\left(\widehat{\mathscr{H}}_{n,m}\right)_{q,c}\to
\operatorname{Pol}\left(\widetilde{X}\right)_q$. Its image will be denoted by
$\operatorname{Pol}\left(\widetilde{\mathscr{H}}_{n,m}\right)_q$. It is easy to verify
that the $*$-algebra $\operatorname{Pol}(\mathscr{H}_{n,m})_q$ is {\it embedded} this way
into $\operatorname{Pol}\left(\widetilde{\mathscr{H}}_{n,m}\right)_q$ and its image is
just the subalgebra in $\operatorname{Pol}\left(\widetilde{\mathscr{H}}_{n,m}\right)_q$
generated by $t_{1j}t_{1k}^*$, $j,k=1,2,\ldots,N$. In what follows we will identify
$\operatorname{Pol}(\mathscr{H}_{n,m})_q$ with its image under the map $J$.

\begin{remark}
\begin{enumerate}
\item $\operatorname{Pol}(\mathscr{H}_{n,m})_q$ can be characterized in
    two ways. Firstly,
$$
\operatorname{Pol}(\mathscr{H}_{n,m})_q=
\left\{\left.f\in\operatorname{Pol}\left(\widetilde{X}\right)_q\right|\:
\triangle_L(f)=1\otimes f\right\}.
$$
Here $\triangle_L$ is the coaction
$\triangle_L:\operatorname{Pol}\left(\widetilde{X}\right)_q\to
\mathbb{C}[\mathfrak{s}(\mathfrak{u}_1\times\mathfrak{u}_{N-1})]_q\otimes
\operatorname{Pol}\left(\widetilde{X}\right)_q$,
$\triangle_L:t_{ij}\mapsto\sum\limits_{k=1}^N\pi(t_{ik})\otimes
t_{kj}$, and $\pi:\operatorname{Pol}\left(\widetilde{X}\right)_q\to
\mathbb{C}[\mathfrak{s}(\mathfrak{u}_1\times\mathfrak{u}_{N-1})]_q$ is
the factorization map with respect to the two-sided ideal in
$\operatorname{Pol}\left(\widetilde{X}\right)_q$ generated by $t_{1k}$,
$t_{k1}$, $k=2,3,\ldots,N$, cf. \cite[11.6.2, 11.6.4]{Kl-Sch}.

\item Another characterization is in observing that
$\operatorname{Pol}(\mathscr{H}_{n,m})_q$ is the subalgebra of
$U_q\mathfrak{s}(\mathfrak{u}_1\times\mathfrak{u}_{N-1})$-invariants under the left
action in $\operatorname{Pol}\left(\widetilde{X}\right)_q$. The latter action is a dual
to the coaction $\triangle_L$ as in \cite[1.3.5, Proposition 15]{Kl-Sch}. To prove the
equivalence one should observe the $U_q
\mathfrak{s}(\mathfrak{u}_1\times\mathfrak{u}_{N-1})$-invariance of $t_{1j}t_{1k}^*$ and
compare the dimensions of graded components of the algebras
$\operatorname{Pol}\left(\widehat{\mathscr{H}}_{n,m}\right)_q$ and $\mathbb
C[GL_N]_q^{U_q \mathfrak{s}(\mathfrak{u}_1\times\mathfrak{u}_{N-1})}$.
\end{enumerate}
\end{remark}

\medskip

We use the notation $t_j$ instead of $t_{1j}$ for the generators of the $*$-algebra
$\operatorname{Pol}\left(\widetilde{\mathscr{H}}_{n,m}\right)_q$.

Let $I_\varphi$, $\varphi\in\mathbb{R}/2\pi\mathbb{Z}$, be the $*$-automorphism of the
$*$-algebra $\operatorname{Pol}\left(\widetilde{\mathscr{H}}_{n,m}\right)_q$ defined on
the generators $\{t_j\}_{j=1,\ldots,N}$ by
\begin{equation}\label{I_phi}
I_\varphi:t_j\mapsto e^{i\varphi}t_j.
\end{equation}

Then one more description of $\operatorname{Pol}(\mathscr{H}_{n,m})_q$ is as follows:
\begin{equation*}\label{polhdef2}
\operatorname{Pol}(\mathscr{H}_{n,m})_q\overset{\mathrm{def}}{=} \left\{\left.f\in
\operatorname{Pol}\left(\widetilde{\mathscr{H}}_{n,m}\right)_q\right|
\:I_\varphi(f)=f\text{ \ for all \ }\varphi\right\}.
\end{equation*}

At the end of this section we list explicit formulas for the action of $U_{q}\mathfrak{su}_{n,m}$ on
$\operatorname{Pol}\left(\widetilde{\mathscr{H}}_{n,m}\right)$.

The action of $U_{q}\mathfrak{su}_{n,m}$ on
$\operatorname{Pol}\left(\widetilde{\mathscr{H}}_{n,m}\right)$ is described as follows:
\begin{equation}\label{act_on_t_i}
\begin{aligned}
E_jt_i &=
\begin{cases}
q^{-1/2}t_{i-1}, & j+1=i,
\\ 0, & \text{otherwise},
\end{cases}&
\\ F_jt_i &=
\begin{cases}
q^{1/2}t_{i+1}, & j=i,
\\ 0, & \text{otherwise},
\end{cases}&
\\ K_j^{\pm 1}t_i &=
\begin{cases}
q^{\pm 1}t_i, & j=i,
\\ q^{\mp 1}t_i, & j+1=i,
\\ t_i, & \text{otherwise},
\end{cases}&
\end{aligned}\hspace{20em}
\end{equation}

\begin{equation}\label{act_on_t_i*}
\begin{aligned}
E_jt_i^* &=
\begin{cases}
-q^{-3/2}t_{i+1}^*, & j=i\;\&\;i\ne n,
\\ q^{-3/2}t_{i+1}^*, & j=i\;\&\;i=n,
\\ 0, & \text{otherwise},
\end{cases}
\\ F_jt_i^* &=
\begin{cases}
-q^{3/2}t_{i-1}^*, & j+1=i\;\&\;i\ne n+1,
\\ q^{3/2}t_{i-1}^*, & j+1=i\;\&\;i=n+1,
\\ 0, & \text{otherwise},
\end{cases}
\\ K_j^{\pm 1}t_i^* &=
\begin{cases}
q^{\mp 1}t_i^*, & j=i,
\\ q^{\pm 1}t_i^*, & j+1=i,
\\ t_i, & \text{otherwise}.
\end{cases}
\end{aligned}\hspace{15em}
\end{equation}

\medskip

\section{\boldmath A $*$-Algebra $\mathscr{D}(\mathscr{H}_{n,m})_q$ of
finite functions}\label{finite}

Let us construct a faithful $*$-representation $T$ of $\operatorname{Pol}(\mathscr{H}_{n,m})_q$ in a
pre-Hilbert space $\mathscr{H}$ (the method of constructing $T$ is well known; see, for example,
\cite{polmat}).

The space $\mathscr{H}$ is a linear span of its orthonormal basis
$\{e(i_1,i_2,\ldots,i_{N-1})|\:i_1,\ldots,i_n\in-\mathbb{Z}_+;\;
i_{n+1},\ldots,i_{N-1}\in\mathbb{N}\}$.

The $*$-representation $T$ is a restriction to $\operatorname{Pol}(\mathscr{H}_{n,m})_q$
of the $*$-representation of
$\operatorname{Pol}\left(\widetilde{\mathscr{H}}_{n,m}\right)$ defined by
\begin{equation}\label{Tj<n}
\begin{aligned}
T(t_j)e(i_1,\ldots,i_{N-1}) &=
q^{\sum\limits_{k=1}^{j-1}i_k}\cdot\left(q^{2(i_j-1)}-1\right)^{1/2}
e(i_1,\ldots,i_j-1,\ldots,i_{N-1}),
\\ T(t_j^*)e(i_1,\ldots,i_{N-1}) &=
q^{\sum\limits_{k=1}^{j-1}i_k}\cdot\left(q^{2i_j}-1\right)^{1/2}
e(i_1,\ldots,i_j+1,\ldots,i_{N-1}),
\end{aligned}
\end{equation}
for $j\le n$,
\begin{equation}\label{Tn<j<N}
\begin{aligned}
T(t_j)e(i_1,\ldots,i_{N-1}) &=
q^{\sum\limits_{k=1}^{j-1}i_k}\cdot\left(1-q^{2(i_j-1)}\right)^{1/2}
e(i_1,\ldots,i_j-1,\ldots,i_{N-1}),
\\ T(t_j^*)e(i_1,\ldots,i_{N-1}) &=
q^{\sum\limits_{k=1}^{j-1}i_k}\cdot\left(1-q^{2i_j}\right)^{1/2}
e(i_1,\ldots,i_j+1,\ldots,i_{N-1}),
\end{aligned}
\end{equation}
for $n<j<N$, and, finally,
\begin{equation}\label{Tj=N}
\begin{aligned}
T(t_N)e(i_1,\ldots,i_{N-1}) &=
q^{\sum\limits_{k=1}^{N-1}i_k}e(i_1,\ldots,i_{N-1}),
\\ T(t_N^*)e(i_1,\ldots,i_{N-1}) &=
q^{\sum\limits_{k=1}^{N-1}i_k}e(i_1,\ldots,i_{N-1}).
\end{aligned}
\end{equation}

Define the elements $\{x_j\}_{j=1,\ldots,N}$ as follows:
\begin{equation}\label{x_j}
x_j\overset{\mathrm{def}}{=}
\begin{cases}
\sum\limits_{k=j}^Nt_kt_k^*, & j>n,
\\ -\sum\limits_{k=j}^nt_kt_k^*+\sum\limits_{k=n+1}^Nt_kt_k^*, & j\le n.
\end{cases}
\end{equation}
Obviously, $x_1=1$, $x_ix_j=x_jx_i$,
\begin{equation}\label{t_jx_k}
t_jx_k=
\begin{cases}
q^2x_kt_j, & j<k,
\\ x_kt_j, & j\ge k,
\end{cases}
\end{equation}
hence
\begin{equation}\label{t_j^*x_k}
t_j^*x_k=
\begin{cases}
q^{-2}x_kt_j^*, & j<k,
\\ x_kt_j^*, & j\ge k.
\end{cases}
\end{equation}

The vectors $e(i_1,\ldots,i_{N-1})$ are joint eigenvectors of the operators
$T(x_j)$, $j=1,2,\ldots,N$:
\begin{equation}\label{Tx_jev}
\begin{aligned}
&T(x_1)=I,&
\\ &T(x_j)e(i_1,\ldots,i_{N-1})=
q^{2\sum\limits_{k=1}^{j-1}i_k}e(i_1,\ldots,i_{N-1}).
\end{aligned}
\end{equation}

The joint spectrum of the pairwise commuting operators $T(x_j)$,
$j=1,2,\ldots,N$, is
\begin{multline*}\label{jsh}
\mathfrak{M}=\left\{(x_1,\ldots,x_N)\in\mathbb{R}^N\right|
\\ \left. x_i/x_j\in q^{2\mathbb{Z}}\;\&\;1=x_1\le x_2\le\ldots\le
x_{n+1}>x_{n+2}>\ldots>x_N>0\right\}.
\end{multline*}

\begin{proposition}\label{faith_pol}
$T$ is a faithful representation of $\operatorname{Pol}(\mathscr{H}_{n,m})_q$.
\end{proposition}

{\bf Proof.} It suffices to verify faithfulness of the (unrestricted) representation $T$
of $\operatorname{Pol}\left(\widetilde{\mathscr{H}}_{n,m}\right)_q$. It is quite obvious
that an arbitrary element of
$\operatorname{Pol}\left(\widetilde{\mathscr{H}}_{n,m}\right)_q$ can be written as a
finite sum
$$
f=\sum_{\fontsize{8}{1}(i_1,\ldots,i_N,j_1,\ldots,j_N):\;i_kj_k=0}
t_1^{i_1}\ldots t_n^{i_n}t_{n+1}^{*i_{n+1}}\ldots
t_N^{*i_N}f_{IJ}(x_2,\ldots,x_N)t_N^{j_N}\ldots
t_{n+1}^{j_{n+1}}t_n^{*j_n}\ldots t_1^{*j_1},
$$
where $f_{IJ}(x_2,\ldots,x_N)$ are polynomials, $I=(i_1,\ldots,i_N)$,
$J=(j_1,\ldots,j_N)$. It follows from the definition of $T$ that every
summand
$$
t_1^{i_1}\ldots t_n^{i_n}t_{n+1}^{*i_{n+1}}\ldots
t_N^{*i_N}f_{IJ}(x_2,\ldots,x_N)t_N^{j_N}\ldots
t_{n+1}^{j_{n+1}}t_n^{*j_n}\ldots t_1^{*j_1}
$$
takes a basis vector $e(k_1,\ldots,k_{N-1})$ to a scalar multiple of the
basis vector
$e(k_1+j_1-i_1,\ldots,k_n+j_n-i_n,k_{n+1}-j_{n+1}+i_{n+1},\ldots,
k_{N-1}-j_{N-1}+i_{N-1})$. Moreover, the sets of indices
$(k_1+j_1-i_1,\ldots,k_{N-1}-j_{N-1}+i_{N-1})$ of the image basis vectors
are different for different monomials, provided the indices of the initial
monomial $e(k_1,\ldots,k_{N-1})$ have modules large enough. Therefore, to
prove our claim, it suffices to choose arbitrarily a summand of $f$ and to
find an initial basis vector $e(k_1,\ldots,k_{N-1})$ in such a way that the
chosen summand does not annihilate (under $T$) the vector
$e(k_1,\ldots,k_{N-1})$.

Let us consider a basis vector $e(k_1,\ldots,k_{N-1})$ with $|k_s|>j_s$ for
all $s=1,\ldots, N-1$. Then
\begin{multline*}
T\left(t_N^{j_N}\ldots t_{n+1}^{j_{n+1}}t_n^{*j_n}\ldots
t_1^{*j_1}\right)e(k_1,\ldots,k_{N-1})=
\\ \mathrm{const}\cdot
e(k_1+j_1,\ldots,k_n+j_n,k_{n+1}-j_{n+1},\ldots,k_{N-1}-j_{N-1}),
\end{multline*}
where $\mathrm{const}\ne 0$.

Moreover, $T(f_{IJ}(x_2,\ldots,x_N))$ acts by multiplying the basis vector
by a (value of a) polynomial $p\left(q^{2k_1},\ldots,q^{2k_{N-1}}\right)$
(due to \eqref{Tx_jev}), where $p(u_1,u_2,\ldots,
u_{N-1})=f_{IJ}(u_1,u_1u_2,\ldots,u_1u_2 \cdots u_{N-1})$, and $p$ is
certainly a nonzero polynomial. A routine argument allows one to find
$k_1,\ldots,k_{N-1}$ such that $|k_s|>j_s$ and
$p\left(q^{2k_1},\ldots,q^{2k_{N-1}}\right)\ne 0$. This proves the claim we
need. \hfill $\square$

\medskip

Let $P$ be the orthogonal projection of $\mathscr{H}$ onto the linear span of vectors
$\{e(\underbrace{0,\ldots,0}_n,i_{n+1},\ldots,i_{N-1})|\:
i_{n+1},\ldots,i_{N-1}\in\mathbb{N}\}$. Of course
$\operatorname{Pol}(\mathscr{H}_{n,m})_q$ does not contain an element $f_0$ such that
$T(f_0)=P$. Our immediate intention is to add $f_0$ with this property.

Consider the $*$-algebra
$\operatorname{Fun}\left(\widetilde{\mathscr{H}}_{n,m}\right)\supset
\operatorname{Pol}\left(\widetilde{\mathscr{H}}_{n,m}\right)$ derived from
$\operatorname{Pol}\left(\widetilde{\mathscr{H}}_{n,m}\right)$ by adding an element $f_0$
to its list of generators and the relations as below to its list of relations:
\begin{equation}\label{f_0rel}
\begin{aligned}
& t_j^*f_0=f_0t_j=0,\qquad j\le n, &
\\ & x_{n+1}f_0=f_0x_{n+1}=f_0, &
\\ & f_0^2=f_0^*=f_0, &
\\ & t_jf_0=f_0t_j;\qquad t_j^*f_0=f_0t_j^*,\qquad j\ge n+1.
\end{aligned}\hspace{10em}
\end{equation}

The relation $I_\varphi f_0=f_0$ allows one to extend the $*$-automorphism $I_\varphi$
\eqref{I_phi} of the algebra
$\operatorname{Pol}\left(\widetilde{\mathscr{H}}_{n,m}\right)$ to the $*$-automorphism of
$\operatorname{Fun}\left(\widetilde{\mathscr{H}}_{n,m}\right)$. Let
\begin{equation*}\label{FunH}
\operatorname{Fun}(\mathscr{H}_{n,m})\overset{\mathrm{def}}{=} \left\{\left.f\in
\operatorname{Fun}\left(\widetilde{\mathscr{H}}_{n,m}\right)\right|\: I_\varphi
f=f\right\}.
\end{equation*}
Obviously, there exists a unique extension of the $*$-representation $T$ to a
$*$-representation of the $*$-algebra $\operatorname{Fun}(\mathscr{H}_{n,m})$ such that
$T(f_0)=P$.

Our subsequent observations involve extensively the two-sided ideal
$\mathscr{D}(\mathscr{H}_{n,m})_q$ of $\operatorname{Fun}(\mathscr{H}_{n,m})$ generated
by $f_0$. We call this ideal the algebra of finite functions on the quantum hyperbolic
space.

\begin{theorem}\label{faith_fin_hyp}
The representation $T$ of $\mathscr{D}(\mathscr{H}_{n,m})_q$ is faithful.
\end{theorem}

{\bf Proof.} Obviously, every $f \in \mathscr{D}(\mathscr{H}_{n,m})_q$ admits a unique
decomposition
\begin{equation*}
f=\sum_{\fontsize{8}{2pt}\begin{array}{l}(i_1\ldots,i_N,j_1\ldots j_N): \\
i_1+\ldots+i_n+j_{n+1}+\ldots+j_N=\\
=j_1+\ldots+j_n+i_{n+1}+\ldots+i_N\end{array}} t_1^{i_1}\ldots
t_n^{i_n}t_{n+1}^{*i_{n+1}}\ldots t_N^{*i_N}f_0t_N^{j_N}\ldots
t_{n+1}^{j_{n+1}}t_n^{*j_n}\ldots t_1^{*j_1}.
\end{equation*}
A straightforward application of the commutation relations \eqref{f_0rel}
allows us to refine the above decomposition as follows:
\begin{equation}\label{fin_ref}
f=\sum_{\fontsize{8}{2pt}\begin{array}{l}(i_1\ldots,i_N,j_1\ldots j_N):\;i_kj_k=0\;\& \\
i_1+\ldots+i_n+j_{n+1}+\ldots+j_N=\\
=j_1+\ldots+j_n+i_{n+1}+\ldots+i_N\end{array}}
t_1^{i_1}\ldots t_n^{i_n}t_{n+1}^{*i_{n+1}}\ldots
t_N^{*i_N} f_{IJ} t_N^{j_N}\ldots
t_{n+1}^{j_{n+1}}t_n^{*j_n}\ldots t_1^{*j_1},
\end{equation}
where
\begin{equation}\label{fin_ref_int}
f_{IJ}=\sum_{K}p_{K}(x_{n+2},\ldots,x_{N-1})t_1^{k_1}t_2^{k_2}\ldots
t_n^{k_n}f_0(t_n^*)^{k_n}\ldots (t_2^*)^{k_2}(t_1^*)^{k_1}
\end{equation}
for some nonzero polynomials $p_{K}$.

Let us consider a basis vector $e(a_1,\ldots,a_{N-1})$. Every summand from
\eqref{fin_ref} takes $e(a_1,\ldots,a_{N-1})$ to a scalar multiple of the
vector $e(a_1+j_1-i_1,\ldots,a_n+j_n-i_n,a_{n+1}-j_{n+1}+i_{n+1},\ldots,
a_{N-1}-j_{N-1}+i_{N-1})$ (nonzero if well defined). By our assumptions on
entries of $I$ and $J$, the subset of nonzero multiples as above are
linearly independent. Thus it suffices to choose arbitrarily a summand in
\eqref{fin_ref} and to prove that it does not annihilate some basis vector.

Let us also choose arbitrarily a summand
$$
p_K(x_{n+2},\ldots,x_{N-1})t_1^{k_1}t_2^{k_2}\ldots
t_n^{k_n}f_0(t_n^*)^{k_n}\ldots (t_2^*)^{k_2}(t_1^*)^{k_1}
$$
from \eqref{fin_ref_int}. Now $T(f_0(t_n^*)^{k_n}\ldots
(t_2^*)^{k_2}(t_1^*)^{k_1})T(t_N^{j_N}\ldots
t_{n+1}^{j_{n+1}}t_n^{*j_n}\ldots
t_1^{*j_1})e(a_1,\ldots,a_{N-1})=\mathrm{const}\cdot
e(0,\ldots,0,a_{n+1}-j_{n+1},\ldots,a_{N-1}-j_{N-1})$. Here $\mathrm{const}
=0$ unless $a_s+k_s+j_s=0$ for $s=1,\ldots,n$ and $a_s>j_s$ for
$s=n+1,\ldots,N-1$. Set $a_s=-k_s-j_s$ for $s=1,\ldots,n$.

Now let us consider the action of $T(p_K(x_{n+2},\ldots,x_{N-1}))$ on
vectors of the form
$e(-k_1,\ldots,-k_n,a_{n+1}-j_{n+1},\ldots,a_{N-1}-j_{N-1})$ with $a_s>j_s$
for $s=n+1,\ldots,N-1$. An argument similar to that used in the final
paragraph of the proof of Proposition \ref{faith_pol} allows us to choose
$a_{n+1},\ldots,a_{N-1}$ in such a way that $T\left(t_1^{i_1}\ldots
t_n^{i_n}t_{n+1}^{*i_{n+1}}\ldots t_N^{*i_N} f_{IJ} t_N^{j_N}\ldots
t_{n+1}^{j_{n+1}}t_n^{*j_n}\ldots t_1^{*j_1}\right)$ does not annihilate
$e(a_1,\ldots,a_{N-1})$. This proves our claim. \hfill $\square$

\begin{remark}
i) Due to \eqref{f_0rel}, $f_0$ can be treated as a function of $x_{n+1}$:
\begin{equation}\label{f_0(x_n+1)}
f_0=f_0(x_{n+1})=
\begin{cases}
1, & x_{n+1}=1,
\\ 0, & x_{n+1}\in q^{-2\mathbb{N}}.
\end{cases}
\end{equation}
(Recall that $\mathrm{spec}\,x_{n+1}=q^{-2\mathbb{Z}_+}$). Thus $f_0$ is a
$q$-analog of the characteristic function of the submanifold
$$
\left\{\left.(t_1,\ldots,t_N)\in\mathbb{C}^N\right|\:t_1=t_2=\ldots=t_n=0
\right\}\cap\mathscr{H}_{n,m}.
$$

ii) Let $f(x_{n+1})$ be a polynomial. Then it follows from \eqref{x_j},
\eqref{t_jx_k} that
\begin{equation}\label{t_ift_^*}
\sum_{i=1}^nt_if(x_{n+1})t_i^*=f(q^2x_{n+1})\sum_{i=1}^nt_it_i^*=
f\left(q^2x_{n+1}\right)(x_{n+1}-1).
\end{equation}
This computation, together with \eqref{f_0(x_n+1)}, allows one to consider
the element $f_1=\sum\limits_{i=1}^nt_if_0t_i^*$ as a function of $x_{n+1}$
such that
$$
f_1(x_{n+1})=
\begin{cases}
q^{-2}-1, & x_{n+1}=q^{-2},
\\ 0, & x_{n+1}=1\text{ \ or \ }x_{n+1}\in q^{-2\mathbb{N}-2}.
\end{cases}
$$
Thus a multiple application of \eqref{t_ift_^*} leads to the following claim:
$\mathscr{D}(\mathscr{H}_{n,m})_q$ contains {\it all} finite functions of $x_{n+1}$ (i.e., such functions $f$
that $f(q^{-n})=0$ for all but finitely many $n\in\mathbb{N}$).

Let us now is to endow $\mathscr{D}(\mathscr{H}_{n,m})_q$ with a structure of
$U_{q}\mathfrak{su}_{n,m}$-module algebra. For that, it suffices to describe the action of the operators
$\left\{E_j,F_j,K_j\right\}_{j=1,\ldots,N-1}$ on $f_0$. Here it is:
\begin{align}
\label{X_n+f_0} E_nf_0 &= -\frac{q^{-1/2}}{q^{-2}-1}t_nf_0t_{n+1}^*, &
\\ \label{X_n-f_0}
F_nf_0 &= -\frac{q^{3/2}}{q^{-2}-1}t_{n+1}f_0t_n^*, &
\\ \label{H_nf_0}
K_nf_0 &= f_0, &
\\ \label{j_ne_n}
E_jf_0 &= F_jf_0=(K_j-1)f_0=0,\qquad j\ne n. &
\end{align}
\end{remark}

\begin{remark}
To see that the above structure of $U_{q}\mathfrak{su}_{n,m}$-module algebra on
$\mathscr{D}(\mathscr{H}_{n,m})_q$ is well-defined, it suffices to use an argument
contained in \cite{polmat}. Here we restrict ourselves to explaining the motives which
lead to \eqref{X_n+f_0} -- \eqref{j_ne_n}. An application of \eqref{act_on_t_i},
\eqref{act_on_t_i*} and \eqref{x_j} allows one to conclude that for any polynomial $f(t)$
\begin{align}\label{X_n^+f}
E_nf(x_{n+1}) &= q^{-1/2}t_n
\frac{f(q^{-2}x_{n+1})-f(x_{n+1})}{q^{-2}x_{n+1}-x_{n+1}}t_{n+1}^*, &
\\ \label{X_n^-f}
F_nf(x_{n+1}) &= q^{3/2}t_{n+1}
\frac{f(q^{-2}x_{n+1})-f(x_{n+1})}{q^{-2}x_{n+1}-x_{n+1}}t_n^*, &
\\ \label{j_ne_n_}
E_jf &= F_jf=(K_j-1)f=0\text{ \ for \ }j\ne n,\;j=1,2,\ldots,N-1.
\end{align}
A subsequent application of \eqref{X_n^+f} -- \eqref{j_ne_n_} to the
non-polynomial function $f_0$ \eqref{f_0(x_n+1)} yields \eqref{X_n+f_0} --
\eqref{j_ne_n}.
\end{remark}

\bigskip

\section{Invariant integral}\label{invint}

The aim of this section is to present an explicit formula for a positive invariant
integral on the space of finite functions $\mathscr{D}(\mathscr{H}_{n,m})_q$ and thereby
to establish its existence.

Let $\nu_q:\mathscr{D}(\mathscr{H}_{n,m})_q\to\mathbb{C}$ be a linear functional defined
by
\begin{equation}\label{nu_q}
\nu_q(f)=\operatorname{Tr}(T(f)\cdot Q)=\int\limits_{\mathscr{H}_{n,m}} fd\nu_q,
\end{equation}
where $Q:\mathscr{H}\to\mathscr{H}$ is the linear operator given on the basis elements
$e(i_1,\ldots,i_{N-1})$ by
\begin{equation}\label{const_int_hyp}
Qe(i_1,\ldots,i_{N-1})= \mathrm{const} \cdot
q^{2\sum\limits_{j=1}^{N-1}(N-j)i_j}e(i_1,\ldots,i_{N-1}), \qquad \mathrm{const}>0.
\end{equation}
Thus $Q=\mathrm{const}\cdot T(x_2\cdot\ldots\cdot x_N)$; this follows from
\eqref{Tx_jev}.

\begin{theorem}
The functional $\nu_q$ determined by \eqref{nu_q} is well defined,
positive, and $U_{q}\mathfrak{su}_{n,m}$-invariant.
\end{theorem}

{\bf Proof.} It follows from \eqref{t_it_j}, \eqref{x_j}, \eqref{t_jx_k} that any element
$f$ of the algebra $\mathscr{D}(\mathscr{H}_{n,m})_q$ can be written in a unique way in
the form
\begin{equation}\label{finite_hyp}
f=\sum_{\fontsize{8}{2pt}\begin{array}{l}(i_1\ldots,i_N,j_1\ldots j_N):\;i_kj_k=0\;\& \\
i_1+\ldots+i_n+j_{n+1}+\ldots+j_N=\\
=j_1+\ldots+j_n+i_{n+1}+\ldots+i_N\end{array}}\mspace{-60mu}
t_1^{i_1}\ldots t_n^{i_n}t_{n+1}^{*i_{n+1}}\ldots
t_N^{*i_N}f_{IJ}(x_2,\ldots,x_N)t_N^{j_N}\ldots
t_{n+1}^{j_{n+1}}t_n^{*j_n}\ldots t_1^{*j_1},
\end{equation}
with $f_{IJ}(x_2,\ldots,x_N)$ being a polynomial in $x_2,\ldots,
x_n,x_{n+2},\ldots,x_N$ and a {\it finite} function in $x_{n+1}$, that is,
$f_{IJ}(x_2,\ldots,x_N)$ has the form
\begin{equation}\label{fIJ}
\sum_\text{finite sum}\alpha_\mathbb{K}x_2^{k_2}\cdots
x_n^{k_n}f_\mathbb{K}(x_{n+1})x_{n+2}^{k_{n+2}}\cdots x_N^{k_N},\qquad
\alpha_\mathbb{K}\in\mathbb{C},
\end{equation}
and $f_\mathbb{K}(q^{-2l})\ne 0$ for finitely many $l\in\mathbb{Z}_+$.

Then, by our definition,
\begin{multline}\label{nu_q_exp}
\nu_q:f\mapsto \mathrm{const} \cdot
\sum_{\fontsize{8}{2pt}\begin{array}{l}(i_1\ldots,i_n)\in(-\mathbb{Z}_+)^n\\
(i_{n+1},\ldots,i_{N-1})\in\mathbb{N}^{m-1}\end{array}}
f_{00}\left(q^{2i_1},q^{2i_1+2i_2},\ldots,q^{2i_1+\ldots+2i_{N-1}}\right)
\cdot
\\ \cdot q^{2(N-1)i_1+2(N-2)i_2+\ldots+2i_{N-1}},
\end{multline}
and the series in the right hand side of \eqref{nu_q_exp} converges for $f$ of the form \eqref{fIJ}.

The positivity of the linear functional $\nu_q$ means that
\begin{equation*}\label{nu_qpos}
\nu_q(f^*f)>0\qquad\text{for}\quad f\ne 0.
\end{equation*}
This follows from the explicit formula \eqref{nu_q_exp} and the {\it faithfulness} of the
$*$-representation $T$ of the algebra $\mathscr{D}(\mathscr{H}_{n,m})_q$ (see Section
\ref{finite}).

What remains is to establish the $U_q\mathfrak{su}_{n,m}$-invariance of
$\nu_q$. The desired invariance is equivalent to
\begin{equation}\label{inv}
\nu_q(E_j f)=0, \qquad \nu_q(F_jf)=0.
\end{equation}
for any $f\in\mathscr{D}(\mathscr{H}_{n,m})_q$ and $j=1,2,\ldots,N-1$. Observe that
$\nu_q$ is a real functional, i.e., $\nu_q(f^*)=\overline{\nu_q(f)}$. The latter relation
follows from selfadjointness of the operator $Q:\mathscr{H}\to\mathscr{H}$ involved in
the definition of $\nu_q$. This allows us to reduce the proof of $\eqref{inv}$ to proving
the abridged version of it
\begin{equation}\label{inv_}
\nu_q(E_jf)=0,\qquad j=1,2,\ldots,N-1.
\end{equation}

We are going to establish \eqref{inv_} for $j\le n$; for other $j$ the
proof is similar.

Moreover, for a function $f$ of the form
$$
f=t_1^{i_1}\ldots t_n^{i_n}t_{n+1}^{*i_{n+1}}\ldots t_N^{*i_N}
f_{IJ}(x_2,\ldots,x_N)t_N^{j_N}\ldots t_{n+1}^{j_{n+1}}t_n^{*j_n}\ldots
t_1^{*j_1}
$$
with $i_kj_k=0$ for $k=1,2,\ldots,N$, one has $\nu_q(E_jf)=0$ if $I\ne\underset{(j+1)\text{th
place}}{(0,\ldots,0,1,0,\ldots,0)}$ and $J\ne\underset{j\text{th place}}{(0,\ldots,0,1,0,\ldots,0)}$ (if
$j<n$) or $I\ne(0,0,\ldots,0)$ and $J\ne\underset{\phantom{\text{places}+1)}j\text{th}\,(j+1)\text{th
places}}{(0,\ldots,0,1,1,0,\ldots,0)}$ (if $j=n$). Thus we have to verify that
$\nu_q\left(E_j\left(t_{j+1}f(x_2,\ldots,x_N)t_j^*\right)\right)=0$.

It can be demonstrated by a direct computation that for $j\le n$
\begin{multline}\label{blf}
E_j(t_{j+1}f(x_2,\ldots,x_N)t_j^*)=
\\ =q^{-1/2}\left[q^2f(x_2,\ldots,x_j,q^2x_{j+1},\ldots,q^2x_N)(x_{j+1}-x_j)
\frac{q^{-2}x_{j+2}-x_{j+1}}{(1-q^2)x_{j+1}}\right.
\\ \left.-f(x_2,\ldots,x_{j+1},q^2x_{j+2},\ldots,q^2x_N)(x_{j+2}-x_{j+1})
\frac{q^{-2}x_{j+1}-x_j}{(1-q^2)x_{j+1}}\right].
\end{multline}

1. Let $j=n$. Then
\begin{multline*}
\nu_q\left(E_j\left(t_{j+1}f(x_2,\ldots,x_N)t_j^*\right)\right)=
\\ =\text{const'}\cdot\mspace{-50mu}\sum_{\fontsize{8}{2pt}
\begin{array}{l}(i_1\ldots,i_n)\in(-\mathbb{Z}_+)^n\\
(i_{n+1},\ldots,i_{N-1})\in\mathbb{N}^{m-1}\end{array}}\mspace{-30mu}
\Bigg[f\left(q^{2i_1},\ldots,q^{2i_1+\ldots+2i_{n-1}},
q^{2i_1+\ldots+2i_n+2},\ldots,q^{2i_1+\ldots+2i_{N-1}+2}\right)\cdot
\\ \cdot\frac{q^2 \left(q^{2i_1+\ldots+2i_n}-q^{2i_1+\ldots+2i_{n-1}}
\right)\left(q^{2i_1+\ldots+2i_{n+1}-2}-q^{2i_1+\ldots+2i_n}\right)}
{q^{2i_1+\ldots+2i_n}}\;-
\\ -f\left(q^{2i_1},\ldots,q^{2i_1+\ldots+2i_n},
q^{2i_1+\ldots+2i_{n+1}+2},\ldots,q^{2i_1+\ldots+2i_{N-1}+2}\right)\cdot
\\ \cdot\frac{\left(q^{2i_1+\ldots+2i_{n+1}}-q^{2i_1+\ldots+2i_n}
\right)\left(q^{2i_1+\ldots+2i_{n}-2}-q^{2i_1+\ldots+2i_{n-1}}\right)}
{q^{2i_1+\ldots+2i_n}}\Bigg]q^{2(N-1)i_1+\ldots+2i_{N-1}}=
\end{multline*}
\begin{multline*}
=\text{const'}\cdot\mspace{-50mu}\sum_{\fontsize{8}{2pt}
\begin{array}{l}(i_1\ldots,i_n)\in(-\mathbb{Z}_+)^n\\
(i_{n+1},\ldots,i_{N-1})\in\mathbb{N}^{m-1}\end{array}}\mspace{-30mu}
\Bigg[f\left(q^{2i_1},\ldots,q^{2i_1+\ldots+2i_{n-1}},
q^{2i_1+\ldots+2i_n+2},\ldots,q^{2i_1+\ldots+2i_{N-1}+2}\right)\cdot
\\ \cdot q^2 \left(q^{2i_n}-1\right)
\left(q^{2i_{n+1}-2}-1\right)\;-
\\ -f\left(q^{2i_1},\ldots,q^{2i_1+\ldots+2i_n},
q^{2i_1+\ldots+2i_{n+1}+2},\ldots,q^{2i_1+\ldots+2i_{N-1}+2}\right)\cdot
\\ \cdot\left(q^{2i_{n+1}}-1 \right)\left(q^{2i_{n}-2}-1\right)
\Bigg]q^{2i_1+\ldots+2i_{n-1}}q^{2(N-1)i_1+\ldots+2i_{N-1}}.
\end{multline*}

Let us consider the inner sum (in $i_n$ and $i_{n+1}$). For brevity, we
denote $f\left(q^{2i_1},\ldots,q^{2i_1+\ldots+2i_{n-1}},
q^{2i_1+\ldots+2i_n+2},\ldots,q^{2i_1+\ldots+2i_{N-1}+2}\right)$ by
$\psi_{i_n+1,i_{n+1}}$.
\begin{multline*}
\mspace{-20mu}\sum\limits_{\fontsize{8}{2pt}
\begin{array}{c}i\in-\mathbb{Z}_+\\
j\in\mathbb{N}\end{array}}\mspace{-20mu} \left[\psi_{i+1,j}\cdot
q^2\left(1-q^{2i}\right)\left(1-q^{2j-2}\right)-\psi_{i,j+1}\cdot
\left(1-q^{2i-2}\right)\left(1-q^{2j}\right)\right]q^{2(N-n)i+2(N-n-1)j}=
\\ =\sum\limits_{i\in-\mathbb{Z}_+,\,j\in\mathbb{N}}
\psi_{i+1,j}\cdot\left(1-q^{2i}\right)\left(1-q^{2j-2}\right)
q^{2(N-n)i+2(N-n-1)j+2}
\\ -\sum\limits_{i\in-\mathbb{Z}_+,\,j\in\mathbb{N}}\psi_{i,j+1}\cdot
\left(1-q^{2i-2}\right)\left(1-q^{2j}\right)q^{2(N-n)i+2(N-n-1)j}
\\ =q^{-2(N-n-1)}\sum_{i\le 1,j\in\mathbb{N}}\psi_{i,j}\left(1-q^{2i-2}\right)
\left(1-q^{2j-2}\right)q^{2(N-n)i+2(N-n-1)j}
\\ -q^{-2(N-n-1)}\sum_{i\in-\mathbb{Z}_+,j\ge 2}\psi_{i,j}
\left(1-q^{2i-2}\right)\left(1-q^{2j-2}\right)q^{2(N-n)i+2(N-n)j}=0.
\end{multline*}
Thus the proof in this case is complete.

2. Let $j<n$.

\begin{multline*}
\mspace{-10mu}\sum_{i,j\in-\mathbb{Z}_+}\mspace{-10mu}
\left[\psi_{i+1,j}\cdot
q^2\left(1-q^{2i}\right)\left(1-q^{2j-2}\right)-\psi_{i,j+1}\cdot
\left(1-q^{2i-2}\right)\left(1-q^{2j}\right)\right]q^{2(N-n)i+2(N-n-1)j}=
\\ =q^{-2(N-n-1)}\sum_{i\le 1,j\in-\mathbb{Z}_+}\psi_{i,j}
\left(1-q^{2i-2}\right)\left(1-q^{2j-2}\right)q^{2(N-n)i+2(N-n-1)j}
\\ - q^{-2(N-n-1)}\sum_{i\in-\mathbb{Z}_+,j\le
1}\psi_{i,j}\left(1-q^{2i-2}\right)\left(1-q^{2j-2}\right)
q^{2(N-n)i+2(N-n)j}=0.
\end{multline*}
The Theorem is proved. \hfill $\square$

\begin{remark}
It is reasonable to choose $\mathrm{const}$ in \eqref{const_int_hyp} so
that the following normalization property is valid:
$$\nu_q(f_0)=1.$$
This allows us to find the constant explicitly:
$$
\mathrm{const}=q^{-(N-n)(N-n-1)} \prod\limits_{j=n+1}^{N-1}\left(1-q^{2(N-j)}\right).
$$
\end{remark}

\bigskip

\section{\boldmath Quantum homogeneous space $\Xi_{n,m}$}

Let $\operatorname{Pol}\left(\widetilde{\Xi}_{n,m}\right)_q$ denotes the quotient algebra of
$\operatorname{Pol}(\widehat{\mathscr{H}}_{n,m})_q$ by the ideal
$\operatorname{Pol}(\widehat{\mathscr{H}}_{n,m})_q\cdot c$ (recall that $c$ belongs to the center of
$\operatorname{Pol}(\widehat{\mathscr{H}}_{n,m})_q$). This is a $q$-analog of the polynomial algebra on the
isotropic cone. Define an automorphism $I_\varphi$, $\varphi\in\mathbb{R}/2\pi\mathbb{Z}$, of the algebra
$\operatorname{Pol}\left(\widetilde{\Xi}_{n,m}\right)_q$ by
\begin{equation*}
I_\varphi(t_j)=e^{i\varphi}t_j,\qquad I_\varphi(t_j^*)=e^{-i\varphi}t_j^*.
\end{equation*}
Then it follows from the definition that
\begin{equation*}\label{polxi}
\operatorname{Pol}(\Xi_{n,m})_q=\left\{\left.f\in
\operatorname{Pol}\left(\widetilde{\Xi}_{n,m}\right)_q\right|\: I_\varphi(f)=f\text{\ \
for any \ }\varphi\right\}.
\end{equation*}
We are going to construct a $*$-representation $T_0$ of the $*$-algebra
$\operatorname{Pol}\left(\widetilde{\Xi}_{n,m}\right)_q$ in a pre-Hilbert space $\mathscr{H}_0$ in such a way
that the restriction of $T_0$ to the subalgebra $\operatorname{Pol}(\Xi_{n,m})_q$ is a faithful
$*$-representation of $\operatorname{Pol}(\Xi_{n,m})_q$.

Let $\{e(i_1,i_2,\ldots,i_{N-1})|\:i_1\in\mathbb{Z};
i_2,\ldots,i_n\in-\mathbb{Z}_+;\; i_{n+1},\ldots,i_{N-1}\in\mathbb{N}\}$ be
the orthonormal basis of the space $\mathscr{H}_0$. Then $T_0$ is defined
as follows.
\begin{equation}\label{T0t1}
\begin{aligned}
T_0(t_1)e(i_1,\ldots,i_{N-1}) &= q^{i_1-1}e(i_1-1,\ldots,i_{N-1}),
\\ T_0(t_1^*)e(i_1,\ldots,i_{N-1}) &= q^{i_1}e(i_1+1,\ldots,i_{N-1}),
\end{aligned}
\end{equation}
\begin{equation}\label{T0j1n}
\left\{
\begin{aligned}
T_0(t_j)e(i_1,\ldots,i_{N-1}) &= q^{\sum\limits_{k=1}^{j-1}i_k}
\left(q^{2(i_j-1)}-1\right)^{1/2}e(i_1,\ldots,i_j-1,\ldots,i_{N-1}),
\\ T_0(t_j^*)e(i_1,\ldots,i_{N-1}) &= q^{\sum\limits_{k=1}^{j-1}i_k}
\left(q^{2i_j}-1\right)^{1/2}e(i_1,\ldots,i_j+1,\ldots,i_{N-1}),
\\ \text{for}\quad 1<j\le n, &
\end{aligned}\right.
\end{equation}
\begin{equation}\label{T0jnN}
\left\{
\begin{aligned}
T_0(t_j)e(i_1,\ldots,i_{N-1}) &= q^{\sum\limits_{k=1}^{j-1}i_k}
\left(1-q^{2(i_j-1)}\right)^{1/2}e(i_1,\ldots,i_j-1,\ldots,i_{N-1}),
\\ T_0(t_j^*)e(i_1,\ldots,i_{N-1}) &= q^{\sum\limits_{k=1}^{j-1}i_k}
\left(1-q^{2i_j}\right)^{1/2}e(i_1,\ldots,i_j+1,\ldots,i_{N-1}),
\\ \text{for}\quad n<j<N,
\end{aligned}\right.
\end{equation}
\begin{equation}\label{T0N}
\begin{aligned}
T_0(t_N)e(i_1,\ldots,i_{N-1}) &= q^{\sum\limits_{k=1}^{N-1}i_k}
e(i_1,\ldots,i_{N-1}),
\\ T_0(t_N^*)e(i_1,\ldots,i_{N-1}) &= q^{\sum\limits_{k=1}^{N-1}i_k}
e(i_1,\ldots,i_{N-1}),
\end{aligned}
\end{equation}

Let us introduce the notation
\begin{equation*}\label{xi_i}
\xi_j=\left\{
\begin{gathered}
\sum_{k=j}^Nt_kt_k^*,\quad j>n,
\\ -\sum_{k=j}^nt_kt_k^*+\sum_{k=n+1}^Nt_kt_k^*,\quad j\le n.
\end{gathered}\right.
\end{equation*}
Evidently, $\xi_1=0$, and the elements $\xi_2,\ldots,\xi_N$ satisfy
\eqref{t_jx_k} -- \eqref{t_j^*x_k} with $x_k$ being replaced by $\xi_k$.
The joint spectrum of the pairwise commuting operators
$\{T_0(\xi_j)\}_{j=\overline{1,N}}$ is the set
\begin{multline*}\label{M0}
\mathfrak{M}_0=\left\{(\xi_1,\ldots,\xi_N)\in\mathbb{R}^N\right|
\\ \left.\xi_j\in q^{2\mathbb{Z}},\;j > 1\;\&\;
0=\xi_1\le\xi_2\le\ldots\le\xi_{n+1}>\xi_{n+2}>\ldots>\xi_N>0\right\}.
\end{multline*}
Similarly to the case of $\operatorname{Pol}(\mathscr{H}_{n,m})_q$, any element from
$\operatorname{Pol}(\Xi_{n,m})_q$ can be written in the form
\begin{equation*}\label{finite_cone}
f=\sum_{\fontsize{8}{2pt}\begin{array}{c}IJ=0 \\ \text{finite sum} \\
i_1+\ldots +i_n+j_{n+1}+\ldots+j_N=\\ =i_{n+1}+\ldots+i_N+j_1+\ldots+j_n \end{array}}
t_1^{i_1}\ldots t_n^{i_n}t_{n+1}^{*i_{n+1}}\ldots
t_N^{*i_N}f_{IJ}(\xi_2,\ldots,\xi_N)t_N^{j_N}\ldots
t_{n+1}^{j_{n+1}}t_n^{*j_n}\ldots t_1^{*j_1},
\end{equation*}
where $f_{IJ}$ are polynomials in $\xi_2,\ldots,\xi_N$, and such
decomposition is unique.

The $*$-algebra $\operatorname{Pol}\left(\widetilde{\Xi}_{n,m}\right)_q$ is a
$U_q\mathfrak{su}_{n,m}$-module algebra. Namely, the action of $U_q\mathfrak{su}_{n,m}$
on the generators $t_j$, $t_j^*$ of
$\operatorname{Pol}\left(\widetilde{\Xi}_{n,m}\right)_q$ is defined by \eqref{act_on_t_i}
-- \eqref{act_on_t_i*}. This definition is correct due to the fact that the element $c$
of the covariant algebra $\operatorname{Pol}(\widetilde{\mathscr{H}}_{n,m})_q$ is
$U_q\mathfrak{su}_{n,m}$-invariant. Thus the $*$-algebra
$\operatorname{Pol}(\Xi_{n,m})_q$ is a $U_q\mathfrak{su}_{n,m}$-module algebra too. The
same computations as in the case of $\operatorname{Pol}(\mathscr{H}_{n,m})_q$ show that
for any polynomial $f(t)$
\begin{equation}\label{uq_acts_on_f(xin+1)}
\begin{aligned}
E_n f(\xi_{n+1}) &= q^{-1/2}t_n
\frac{f(q^{-2}\xi_{n+1})-f(\xi_{n+1})}{q^{-2}\xi_{n+1}-\xi_{n+1}}t_{n+1}^*,
\\ F_n f(\xi_{n+1}) &= q^{3/2}t_{n+1}
\frac{f(q^{-2}\xi_{n+1})-f(\xi_{n+1})}{q^{-2}\xi_{n+1}-\xi_{n+1}}t_n^*,
\\ (K_n-1)f(\xi_{n+1}) &= E_j f(\xi_{n+1})=F_j f(\xi_{n+1})=(K_j-1)f(\xi_{n+1})=0,\quad j\ne n.
\end{aligned}
\end{equation}

Now \eqref{t_jx_k}, \eqref{t_j^*x_k}, and \eqref{uq_acts_on_f(xin+1)} allow one to
introduce the covariant $*$-algebra $\mathscr{D}(\Xi_{n,m})$ of finite functions on the
quantum homogeneous space $\Xi_{n,m}$. It is formed by elements of the form
\eqref{finite_hyp} with $\xi_k$ instead of $x_k$, where $f_{IJ}(\xi_2,\ldots,\xi_N)$ are
polynomials of $\xi_2,\ldots,\xi_n,\xi_{n+2},\ldots,\xi_N$ and finite functions of
$\xi_{n+1}$ (i.e., $f_{IJ}$ has the form \eqref{fIJ} where $f_\mathbb{K}(q^{2l})\ne 0$
for finitely many $l\in\mathbb{Z}$).

\begin{theorem}\label{faith_fin_cone}
$T_0$ can be extended to a faithful $*$-representation of the $*$-algebra
$\mathscr{D}(\Xi_{n,m})$.
\end{theorem}

\begin{remark}
The algebra $\operatorname{Pol}(\mathscr{H}_{n,m})_q$ has the same list of generators as
$\operatorname{Pol}(\widetilde{\Xi})_q$ while the lists of relations differ by replacing
$c-1=0$ with $c=0$. Furthermore, the differences between the formulas \eqref{Tj<n} --
\eqref{Tj=N} and \eqref{T0t1} -- \eqref{T0N} are low enough to enable us to apply the
same argument in proving Theorems \ref{faith_fin_cone} and \ref{faith_fin_hyp}.
\end{remark}

Now let us construct an invariant integral on $\mathscr{D}(\Xi_{n,m})$. Denote by $\nu_q^0$ the linear
functional $\nu_q^0:\mathscr{D}(\Xi_{n,m})\to\mathbb{C}$ given by
\begin{equation}\label{n_q0}
\nu_q^0(f)=\operatorname{Tr}(T_0(f)\cdot Q_0)\;
\left(=\int\limits_{\Xi_{n,m}}fd\nu_q^0\right)
\end{equation}
with $Q_0:\mathscr{H}_0\to\mathscr{H}_0$ being the linear map given by
\begin{equation}\label{const_int_con}
Q_0e(i_1,\ldots,i_{N-1})= \mathrm{const}\cdot
q^{2\sum\limits_{j=1}^{N-1}(N-j)i_j}e(i_1,\ldots,i_{N-1}).
\end{equation}

\begin{theorem}\label{nuq0_inv}
The functional $\nu_q^0$ is well-defined, positive, and
$U_q\mathfrak{su}_{n,m}$-invariant.
\end{theorem}

{\bf Proof.} It follows from the definition that
\begin{equation}\label{nu_q^0_exp}
\nu_q^0(f)=\mathrm{const}\cdot\mspace{-50mu}\sum_{\fontsize{8}{2pt}
\begin{array}{c}i_1\in\mathbb{Z}\\ (i_2\ldots,i_n)\in(-\mathbb{Z}_+)^{n-1}\\
(i_{n+1},\ldots,i_{N-1})\in\mathbb{N}^{m-1}\end{array}}\mspace{-50mu}
f_{00}\left(q^{2i_1},q^{2i_1+2i_2},\ldots,q^{2i_1+\ldots+2i_{N-1}}\right)
q^{2i_1(N-1)+\ldots+2i_{N-1}}.
\end{equation}
Here $f_{00}$ is the function involved in the decomposition
\eqref{finite_hyp} of $f$.

To prove that the definition \eqref{n_q0} of $\nu_q^0$ is correct, it now
suffices to show that the series in the r.h.s. of \eqref{nu_q^0_exp} is
absolutely convergent for $f_{00}$ satisfying the condition
\begin{equation*}
f_{00}\left(\xi_2,\ldots,\xi_n,q^{2l},\xi_{n+2},\ldots,\xi_N\right)=0\quad
\text{for}\;l\ne l_0.
\end{equation*}
Let $f_{00}$ be such a function. Then
\begin{multline}\label{series}
\sum_{\fontsize{8}{2pt}\begin{array}{c}i_1\in\mathbb{Z}\\
(i_2\ldots,i_n)\in(-\mathbb{Z}_+)^{n-1}\\
(i_{n+1},\ldots,i_{N-1})\in\mathbb{N}^{m-1}\end{array}}
f_{00}\left(q^{2i_1},q^{2i_1+2i_2},\ldots,q^{2i_1+\ldots+2i_n}
q^{2i_1+\ldots+2i_{N-1}}\right)q^{2i_1(N-1)+\ldots+2i_{N-1}}=
\\ =\sum_{\fontsize{8}{2pt}\begin{array}{c}
(i_2\ldots,i_n)\in(-\mathbb{Z}_+)^{n-1}\\
(i_{n+1},\ldots,i_{N-1})\in\mathbb{N}^{m-1}\end{array}}
f_{00}\left(q^{2l_0-2i_2-\ldots-2i_n},q^{2l_0-2i_3-\ldots-2i_n},\ldots,
q^{2l_0-2i_n},q^{2l_0},q^{2l_0+2i_{n+1}},\ldots\right)\cdot
\\ \cdot q^{2l_0(N-1)}\cdot q^{2i_1(N-1)+\ldots+2i_{N-1}}\cdot
q^{-2i_2-4i_3-\ldots-2(n-1)i_n}\cdot
q^{2i_{n+1}(m-1)+i_{n+2}(m-2)+\ldots+2i_{N-1}}.
\end{multline}
It is implicit here that only terms with $i_1+\ldots+i_n=l_0$ can be
non-zero; also, the following obvious equality is used:
$$
q^{2(N-1)i_1+\ldots+2i_{N-1}}=q^{2i_1}\cdot q^{2i_1+2i_2}\cdot\ldots\cdot
q^{2i_1+\ldots+2i_{N-1}}.
$$

Now to establish the convergence of the series \eqref{series}, it suffices
to recall that $f_{00}$ is a polynomial in
$\xi_2,\ldots,\xi_n,\xi_{n+2},\ldots,\xi_N$.

The positive definiteness of $\nu_q^0$ can be explained in the same way as it was done in Section
\ref{invint} for $\nu_q$.

Let us turn to proving the invariance of $\nu_q^0$. To do this, one needs
to reproduce the proof of a similar fact for $\nu_q$ almost literally,
including the computations of cases 1 and 2. But now there is one more case
to be considered:

3. Let $j=1$, then (see \eqref{blf})
\begin{multline*}
E_1(t_2f(\xi_2,\ldots,\xi_N)t_1^*)=
\\ =q^{-1/2}\left[f(q^2\xi_2,\ldots,q^2\xi_N)
\frac{\xi_2(\xi_3-q^2\xi_2)}{(1-q^2)\xi_2} - f(\xi_2,q^2\xi_3,\ldots,q^2\xi_N)
\frac{q^{-2}\xi_2(\xi_3-\xi_2)}{(1-q^2)\xi_2}\right]=
\\ =\frac{q^{-1/2}}{1-q^2}\left[f(q^2\xi_2,\ldots,q^2\xi_N)
(\xi_3-q^2\xi_2)-q^{-2}f(\xi_2,q^2\xi_3,\ldots,q^2\xi_N)(\xi_3-\xi_2)\right].
\end{multline*}

Now let us show that $\nu_q^0(E_1(t_2f(\xi_2,\ldots,\xi_N)t_1^*))=0$. In fact,
\begin{multline}\label{nuX()=0}
\nu_q^0(E_1(t_2f(\xi_2,\ldots,\xi_N)t_1^*))=
\\ =\mathrm{const}'\cdot
\sum_{\fontsize{8}{2pt}\begin{array}{c}i_1\in\mathbb{Z}\\
(i_2\ldots,i_n)\in(-\mathbb{Z}_+)^{n-1}\\
(i_{n+1},\ldots,i_{N-1})\in\mathbb{N}^{m-1}\end{array}}
\left[f(q^{2i_1+2},q^{2i_1+2i_2+2},\ldots,q^{2i_1+\ldots+2i_{N-1}+2})
(q^{2i_2-2}-1)q^{2i_1+2}\right.
\\ \left.-f(q^{2i_1},q^{2i_1+2i_1+2},\ldots,q^{2i_1+\ldots+2i_{N-1}+2})q^{-2}
(q^{2i_2}-1)q^{2i_1}\right]q^{2i_1(N-1)+\ldots+2i_{N-1}}.
\end{multline}
As usual, we denote $f(q^{2i_1+2},q^{2i_1+2i_2+2},\ldots,q^{2i_1+\ldots+2i_{N-1}+2})$ by
$\psi_{i_1+1,i_2}$. Let us compute the inner sum over $i_1$ and $i_2$ in the r.h.s. of
\eqref{nuX()=0}.
\begin{multline*}
\sum_{i\in\mathbb{Z}, j\in-\mathbb{Z}_+} \left[q^2\psi_{i+1,j}(q^{2j-2}-1)-
q^{-2}\psi_{i,j+1}(q^{2j}-1)\right]\cdot q^{2iN} q^{2j(N-2)}=
\\ =\sum_{i\in\mathbb{Z}, j\in-\mathbb{Z}_+} \psi_{i,j}(q^{2j-2}-1)\cdot
q^{2iN+2jN-4j-2N+2}-\sum_{i\in\mathbb{Z}, j \leq 1} \psi_{i,j}(q^{2j-2}-1)\cdot
q^{2iN+2j(N-2)-2N+2}=0. \;\square
\end{multline*}

\begin{remark}
Here $\mathrm{const}$ is chosen in \eqref{const_int_con} so that the
following normalization property is valid:
$$\nu_q^0(f_0)=1.$$
This allows us to find the constant explicitly:
$$
\mathrm{const}=q^{-(N-n)(N-n-1)}\prod\limits_{j=1}^{n-1}\left(1-q^{2j}\right)
\prod\limits_{j=1}^{N-n-1}\left(1-q^{2j}\right).
$$
\end{remark}

\bigskip

\section{Principal non-unitary and unitary series of representations of
\boldmath $U_q\mathfrak{su}_{n,m}$ related to the space $\Xi_{n,m}$}

The element $\xi_{n+1}$ quasi-commutes with all the generators of the algebra
$\operatorname{Pol}(\Xi_{n,m})_q$. Thus $(\xi_{n+1})^{\mathbb Z_+}$ is an Ore set and one
can consider a localization $\operatorname{Pol}(\Xi_{n,m})_{q,\xi_{n+1}}$ of the algebra
$\operatorname{Pol}(\Xi_{n,m})_q$ with respect to the multiplicative system
$(\xi_{n+1})^{\mathbb Z_+}$. Evidently, the $U_q \mathfrak{su}_{n,m}$-module algebra
structure extends to the localization in a unique way.

Denote by $\gamma$ the automorphism of the algebra
$\operatorname{Pol}\left(\widetilde{\Xi}_{n,m}\right)_q$ given on the generators by
\begin{equation*}
\gamma:t_j\mapsto qt_j,\qquad t_j^*\mapsto qt_j^*.
\end{equation*}
Note that $\gamma$ is well defined due to the homogeneity of the defining relations for
$\operatorname{Pol}\left(\widetilde{\Xi}_{n,m}\right)_q$. Obviously,
$\gamma(\xi_{n+1})=q^2\xi_{n+1}$, and this allows one to extend $\gamma$ to an
automorphism of the algebra $\operatorname{Pol}(\Xi_{n,m})_{q,\xi_{n+1}}$, which commutes
with the action of $U_q\mathfrak{su}_{n,m}$. This can be deduced from \eqref{act_on_t_i},
\eqref{act_on_t_i*}, and \eqref{uq_acts_on_f(xin+1)}.



Introduce the $*$-algebra $\mathscr{E}(\Xi_{n,m})_q$ of elements of the form
\begin{equation*}
f=\sum_{\fontsize{8}{2pt}\begin{array}{c}IJ=0 \\ \text{finite sum} \\
i_1+\ldots +i_n+j_{n+1}+\ldots+j_N=\\ =i_{n+1}+\ldots+i_N+j_1+\ldots+j_n \end{array}} t_1^{i_1}\ldots
t_n^{i_n}t_{n+1}^{*i_{n+1}}\ldots t_N^{*i_N}f_{IJ}(\xi_2,\ldots,\xi_N)t_N^{j_N}\ldots
t_{n+1}^{j_{n+1}}t_n^{*j_n}\ldots t_1^{*j_1},
\end{equation*}
with
\begin{equation}\label{fIJ_cone}
f_{IJ}(\xi_2,\ldots,\xi_N)=\sum_{\fontsize{8}{2pt}
\begin{array}{c}\text{finite sum}\\ k_2, \ldots, k_n, k_{n+2},\ldots,k_N
\in \mathbb{Z}_+\\ k_{n+1}\in\mathbb{C}
\end{array}}\alpha_\mathbb{K}\xi_2^{k_2}\xi_3^{k_3}\ldots\xi_N^{k_N}.
\end{equation}
Here $\alpha_\mathbb{K}\in\mathbb{C}$ and the algebra structure is given by \eqref{t_jx_k}, \eqref{t_j^*x_k}.

Given $s\in\mathbb{C}$, let $\mathscr{E}_s(\Xi_{n,m})_q$ be the subspace in $\mathscr{E}(\Xi_{n,m})_q$ of
those elements which have the `homogeneity degree' equal to $s-N+1$:
\begin{equation}\label{deg_s-N+1}
\gamma(f)=q^{s-N+1}\cdot f.
\end{equation}
Thus $\mathscr{E}_s(\Xi_{n,m})_q$ is a $U_q\mathfrak{su}_{n,m}$-submodule in $\mathscr{E}(\Xi_{n,m})_q$. We
call these submodules the modules of the principal non-unitary series related to $\Xi_{n,m}$.

Now let us construct an invariant integral in $\mathscr{E}_{-N+1}(\Xi_{n,m})_q$.

Note that $\mathscr{D}(\Xi_{n,m})_q$ can be made a covariant
$\mathscr{E}(\Xi_{n,m})_q$-bimodule using the relations \eqref{t_jx_k}, \eqref{t_j^*x_k}.

Let $\chi_l\in\mathscr{D}(\Xi_{n,m})_q$ be the function of $\xi_{n+1}$ such that
\begin{equation*}\label{chi_l}
\chi_l(q^{2k})=\delta_{kl},\qquad k,l\in\mathbb{Z}.
\end{equation*}

\begin{lemma}\label{bq}
For any $f\in\mathscr{E}_{-N+1}(\Xi_{n,m})_q$, the integral
\begin{equation}\label{bql}
b_q^{(l)}(f)\stackrel{\mathrm{def}}{=} \int\limits_{\Xi_{n,m}}f\cdot\chi_ld\nu_q^0
\end{equation}
does not depend on $l$.
\end{lemma}

{\bf Proof.}
\begin{multline}\label{bql_exp}
b_q^{(l)}(f)=
\\ =\mathrm{const}\mspace{-50mu}\sum_{\fontsize{8}{2pt}
\begin{array}{c}i_1\in\mathbb{Z}\\ (i_2\ldots,i_n)\in(-\mathbb{Z}_+)^{n-1}
\\ (i_{n+1},\ldots,i_{N-1})\in\mathbb{N}^{m-1}\end{array}}\mspace{-50mu}
f_{00}(q^{2i_1},q^{2i_1+2i_2},\ldots,q^{2i_1+\ldots+2i_{N-1}})
\chi_l(q^{2i_1+\ldots+2i_{N-1}})q^{2i_1(N-1)+\ldots+2i_{N-1}}=
\\ =\mathrm{const}\mspace{-30mu}\sum_{\fontsize{8}{2pt}
\begin{array}{c}(i_2\ldots,i_n)\in(-\mathbb{Z}_+)^{n-1}
\\ (i_{n+1},\ldots,i_{N-1})\in\mathbb{N}^{m-1}\end{array}}\mspace{-30mu}
f_{00}(q^{2l-2i_2-\ldots-2i_n},q^{2l-2i_3-\ldots-2i_n},\ldots,q^{2l-2i_n},
q^{2l},q^{2l+2i_{n+1}},\ldots)\cdot \\ \cdot q^{2l(N-1)}\cdot
q^{-2i_2-4i_3-\ldots-2(n-1)i_n+2i_{n+1}(m-1)+2i_{n+2}(m-2)+\ldots+2i_{N-1}}.
\end{multline}
Clearly, $f\in\mathscr{E}_{-N+1}(\Xi_{n,m})_q$ implies
$$\gamma(f_{00}(\xi_2,\ldots,\xi_N))=q^{-2N+2}f_{00}(\xi_2,\ldots,\xi_N),$$
or, equivalently,
$$f_{00}(q^2\xi_2,\ldots,q^2\xi_N)=q^{-2N+2}f_{00}(\xi_2,\ldots,\xi_N),$$
and thus the r.h.s. of \eqref{bql_exp} can be rewritten as follows
\begin{multline}\label{bql_exp_}
\mathrm{const}\mspace{-20mu}\sum_{\fontsize{8}{2pt}
\begin{array}{c}(i_2\ldots,i_n)\in(-\mathbb{Z}_+)^{n-1}
\\ (i_{n+1},\ldots,i_{N-1})\in\mathbb{N}^{m-1}\end{array}}\mspace{-20mu}
q^{2l(N-1)}f_{00}(q^{-2i_2-\ldots-2i_n},q^{-2i_3-\ldots-2i_n},\ldots,
q^{-2i_n},1,q^{2i_{n+1}},\ldots)\cdot
\\ \cdot q^{2l(N-1)}\cdot
q^{-2i_2-4i_3-\ldots-2(n-1)i_n+2i_{n+1}(m-1)+2i_{n+2}(m-2)+\ldots+2i_{N-1}}=
\\ =\mathrm{const}\mspace{-20mu}\sum_{\fontsize{8}{2pt}
\begin{array}{c}(i_2\ldots,i_n)\in(-\mathbb{Z}_+)^{n-1}
\\ (i_{n+1},\ldots,i_{N-1})\in\mathbb{N}^{m-1}\end{array}}\mspace{-20mu}
f_{00}(q^{-2i_2-\ldots-2i_n},q^{-2i_3-\ldots-2i_n},\ldots,q^{-2i_n},1,
q^{2i_{n+1}},\ldots)\cdot
\\ \cdot
q^{-2i_2-4i_3-\ldots-2(n-1)i_n+2i_{n+1}(m-1)+2i_{n+2}(m-2)+\ldots+2i_{N-1}}.
\;\square
\end{multline}

Introduce the notation $b_q(f)$ or $\int fdb_q$ for the linear functional \eqref{bql} on
$\mathscr{E}_{-N+1}(\Xi_{n,m})_q$. It follows from the proof of Lemma \ref{bq} that
\begin{multline}\label{bq_exp}
b_q(f)=(q^{-2}-1)^N\cdot
\\ \cdot\mspace{-30mu}\sum_{\fontsize{8}{2pt}
\begin{array}{c}(j_1\ldots,j_{n-1})\in(-\mathbb{Z}_+)^{n-1}
\\ (i_1,\ldots,i_{m-1})\in\mathbb{N}^{m-1}\end{array}}\mspace{-40mu}
f_{00}(q^{2j_1+\ldots+2j_{n-1}},q^{2i_2+\ldots+2j_{n-1}},\ldots,q^{-2j_{n-1}},
1,q^{2i_1},q^{2i_1+2i_2},\ldots,q^{2i_1+\ldots,2i_{m-1}})\cdot
\\ \cdot
q^{2j_1+4j_2+\ldots+2(n-1)j_{n-1}}\cdot
q^{2(m-1)i_1+2(m-2)i_2+\ldots+2i_{m-1}}.
\end{multline}

\begin{theorem}
$b_q$ is an invariant integral on $\mathscr{E}_{-N+1}(\Xi_{n,m})_q$.
\end{theorem}

{\bf Proof.} By \eqref{uq_acts_on_f(xin+1)}, the functions of $\xi_{n+1}$ are
$U_q\mathfrak{s}(\mathfrak{u}_n\times\mathfrak{u}_m)$-invariants. Thus $b_q$ is a
$U_q\mathfrak{s}(\mathfrak{u}_n\times\mathfrak{u}_m)$-invariant functional (see Theorem
\ref{nuq0_inv}). It remains to prove that $b_q\left(F_n f\right)=b_q\left(E_n f\right)=0$
for $f\in\mathscr{E}_{-N+1}(\Xi_{n,m})_q$. Let us prove just one of these two equalities,
for example, $b_q(E_nf)=\int\limits_{\Xi_{n,m}}E_nf\cdot\chi_ld\nu_q^0=0$.

The invariance of $\nu_q^0$ and the fact that $\mathscr{D}(\Xi_{n,m})_q$ is a covariant
$\mathscr{E}(\Xi_{n,m})_q$-bimodule imply that
$$b_q(E_nf)=-q^{-1}\int f\cdot E_n\chi_ld\nu_q^0, \qquad f \in
\mathscr{E}_{-N+1}(\Xi_{n,m})_q$$ (the integration by parts is used here, see
\cite[Chapter 4]{Ch-P}).

By \eqref{uq_acts_on_f(xin+1)},

\begin{multline}\label{lat_expr}
-q^{-1}\int f\cdot E_n\chi_ld\nu_q^0=-q^{-1}\int f\cdot q^{-1/2}t_n
\frac{\chi_l(q^{-2}\xi_{n+1})-\chi_l(\xi_{n+1})}{(q^{-2}-1)\xi_{n+1}}
t_{n+1}^*d\nu_q^0=
\\ =-\frac{q^{-3/2}}{(q^{-2}-1)}\int f\cdot t_n
\frac{\chi_{l+1}(\xi_{n+1})-\chi_l(\xi_{n+1})}{\xi_{n+1}}t_{n+1}^*d\nu_q^0=
\\ =-\frac{q^{-3/2}}{(q^{-2}-1)}\operatorname{Tr}\left[T_0\left(f\cdot t_n
\frac{\chi_{l+1}-\chi_l}{\xi_{n+1}}t_{n+1}^*\right)Q_0\right]=
\\ =-\frac{q^{-3/2}}{(q^{-2}-1)}(q^{-2}-1)^N\operatorname{Tr}
\left[T_0\left(f\cdot
t_n\frac{\chi_{l+1}-\chi_l}{\xi_{n+1}}t_{n+1}^*\xi_2\xi_3\ldots\xi_N\right)
\right]=
\\ =\mathrm{const}(q,n,N)\operatorname{Tr}\left[T_0\left(f\cdot
t_n\frac{\chi_{l+1}-\chi_l}{\xi_{n+1}}\xi_2\xi_3\ldots\xi_Nt_{n+1}^*\right)
\right]=
\\ =\mathrm{const}(q,n,N)\operatorname{Tr}\left[T_0\left(t_{n+1}^*f\cdot
t_n\frac1{\xi_{n+1}}(\chi_{l+1}-\chi_l)\xi_2\xi_3\ldots\xi_N\right)\right]=
\\ =\mathrm{const}'(q,n,N)\operatorname{Tr}\left[T_0\left(t_{n+1}^*f\cdot
t_n\frac1{\xi_{n+1}}(\chi_{l+1}-\chi_l)Q_0\right)\right]=
\\ =\mathrm{const}'(q,n,N)\int t_{n+1}^*f\cdot
t_n\frac1{\xi_{n+1}}(\chi_{l+1}-\chi_l)d\nu_1^0.
\end{multline}

If $f\in\mathscr{E}_{-N+1}(\Xi_{n,m})_q$, one has $t_{n+1}^*f\cdot
t_n\frac1{\xi_{n+1}}\in\mathscr{E}_{-N+1}(\Xi_{n,m})_q$. Thus the latter expression in
\eqref{lat_expr} can be rewritten as follows:
\begin{multline*}
\mathrm{const}'(q,n,N)\left(\int t_{n+1}^*f\cdot
t_n\frac1{\xi_{n+1}}\chi_{l+1}d\nu_1^0-\int t_{n+1}^*f\cdot
t_n\frac1{\xi_{n+1}}\chi_ld\nu_1^0\right)=
\\ =\mathrm{const}'(q,n,N)\left(b_q^{(l+1)}\left(t_{n+1}^*f\cdot
t_n\frac1{\xi_{n+1}}\right)-b_q^l\left(t_{n+1}^*f\cdot
t_n\frac1{\xi_{n+1}}\right)\right).
\end{multline*}

It follows from Lemma \ref{bq} that the latter difference is zero. \hfill
$\square$

If $f_1\in\mathscr{E}_s(\Xi_{n,m})_q$ and $f_2\in\mathscr{E}_{-s}(\Xi_{n,m})_q$, one has
$f_1\cdot f_2\in\mathscr{E}_{-N+1}(\Xi_{n,m})_q$. Now an application of the standard
arguments (see, e.g., \cite[Chapter 4]{Ch-P}) which set correspondence between invariant
integrals and invariant pairings, yields

\begin{corollary}\label{pair_inv}
The pairing $\mathscr{E}_s(\Xi_{n,m})_q\times
\mathscr{E}_{-s}(\Xi_{n,m})_q\to\mathbb{C}$,
$$
(f_1,f_2)\mapsto\langle f_1,f_2\rangle\stackrel{\mathrm{def}}{=}\int
f_1f_2db_q
$$
is $U_q\mathfrak{su}_{n,m}$-invariant.
\end{corollary}

Obviously, the involution $*$ of the $*$-algebra $\mathscr{E}(\Xi_{n,m})_q$ maps
$\mathscr{E}_{i\lambda}(\Xi_{n,m})_q$ to $\mathscr{E}_{-i\lambda}(\Xi_{n,m})_q$ for
$\lambda\in\mathbb{R}$.

\begin{proposition}
The sesquilinear form
\begin{equation}\label{seql_f}
(f_1,f_2)=\int f_2^*f_1db_q,\qquad f_1,f_2\in\mathscr{E}_{i\lambda}(\Xi_{n,m})_q,
\end{equation}
is invariant and positive definite.
\end{proposition}

{\bf Proof.} The invariance follows immediately from Corollary \ref{pair_inv} (the
standard arguments from \cite[Chapter 4]{Ch-P} are to be applied here again).

To see that the form \eqref{seql_f} is positive definite, one should recall that the
integral $\nu_q^0$ is positive definite (Theorem \ref{nuq0_inv}), and use the following
computations:
\begin{multline*}
(f,f)=\int f^*fdb_q=\int\limits_{\Xi_{n,m}}f^*f\chi_ld\nu_q^0=
\operatorname{Tr}\left(T_0\left(f^*f\chi_l\right)Q_0\right)=
\operatorname{Tr}\left(T_0\left(f^*f\chi_l\chi_l\right)Q_0\right)=
\\ =\operatorname{Tr}\left(T_0\left(f^*f\chi_l\cdot\mathrm{const}\cdot\xi_2
\ldots\xi_N\chi_l\right)\right)=
\operatorname{Tr}\left(T_0\left(\chi_lf^*f\chi_l\right)Q_0\right)=
\operatorname{Tr}\left(T_0\left(\chi_l^*f^*f\chi_l\right)Q_0\right)=
\\ =\int\limits_{\Xi_{n,m}}(f\chi_l)^*f\chi_ld\nu_q^0.
\end{multline*}
Here $f\in\mathscr{E}_{i\lambda}(\Xi_{n,m})_q$, $\lambda\in\mathbb{Z}$, and the obvious relations
$\chi_l^2=\chi_l$, $\chi_l^*=\chi_l$, $\chi_l\xi_k=\xi_k\chi_l$ are used. \hfill $\square$

Thus $\mathscr{E}_{i\lambda}(\Xi_{n,m})_q$, $\lambda\in\mathbb{R}$, are unitary
$U_q\mathfrak{su}_{n,m}$-modules. They will be called the modules of the principal
unitary series related to $\Xi_{n,m}$.

\bigskip

Let us look at the structure of $\mathscr{E}_\lambda(\Xi_{n,m})_q$ as a
$U_q(\mathfrak{sl}_n\times\mathfrak{sl}_m)$-module. Let $L^{(n)}(\lambda)$ be the finite dimensional simple
$U_q \mathfrak{sl}_n$-module with highest weight $\lambda$. Also let $\varpi_j, j=1,\ldots,n-1,$ be the
fundamental weights of the Lie algebra $\mathfrak{sl}_n$.

Now we recall that if $A$ is a Hopf algebra and $V_1$ is an $A$-module, and $B$ is a Hopf algebra and $V_2$
is a $B$-module then $V_1 \boxtimes V_2$ denotes their tensor product endowed with the structure of $A\otimes
B$-module in the natural way.

\begin{theorem}
The $U_q(\mathfrak{sl}_n\times\mathfrak{sl}_m)$-module $\mathscr{E}_{2s}(\Xi_{n,m})_q$ splits as a
multiplicity free direct sum of its simple submodules
\begin{equation}\label{subm}
L^{(n)}(k\omega_1+l\omega_{n-1})\boxtimes L^{(m)}(l'\omega_1+k'\omega_{m-1}),
\end{equation}
with $k,l,k',l'\ge 0$, $k+l'=k'+l$. Every such submodule is generated by the highest weight vector of the form
\begin{equation}\label{uqslN_highest}
t_1^kt_N^{*k'}\xi_{n+1}^{(s-k'-l)}t_{n+1}^{l'}t_n^{*l}.
\end{equation}
\end{theorem}
{\bf Proof.} For simplicity, we prove the claim in the special case $s=(N-1)/2$, the other cases are similar.
Each element $f \in \mathscr{E}_{N-1}(\Xi_{n,m})_q$ can be decomposed in the following way
\begin{equation*}
f= \sum \limits_{\substack{\text{finite sum} \\
i_1+\ldots +i_n+j_{n+1}+\ldots+j_N=\\
=i_{n+1}+\ldots+i_N+j_1+\ldots+j_n=\lambda}}
\!\!\!\!\!\!\!\!\!\!t_1^{i_1}\ldots t_n^{i_n}t_{n+1}^{*i_{n+1}}\ldots
t_N^{*i_N} \cdot \xi_{n+1}^{-\lambda} \cdot t_N^{j_N}\ldots
t_{n+1}^{j_{n+1}}t_n^{*j_n}\ldots t_1^{*j_1}.
\end{equation*}
Evidently, in all such decompositions $\lambda \in\mathbb Z_+$. For every
fixed decomposition of $f$ let us consider the largest $\lambda$ through
all the terms, and then denote by $\lambda(f)$ the smallest $\lambda$
throughout all such decompositions of $f$. Now we introduce a filtration
\begin{equation*}
\mathscr{E}_{N-1}(\Xi_{n,m})_q=\bigcup_{a=0}^\infty
\mathscr{E}_{N-1}(\Xi_{n,m})_{q,a},
\end{equation*}
where
\begin{equation*}
\mathscr{E}_{N-1}(\Xi_{n,m})_{q,a}=\{f \in \mathscr{E}_{N-1}(\Xi_{n,m})_q
\,\,| \,\, \lambda(f) \leq a\}.
\end{equation*}

By obvious reasons,
$t_1^kt_N^{*k'}\xi_{n+1}^{((N-1)/2-k'-l)}t_{n+1}^{l'}t_n^{*l}$ generates a
$U_q \mathfrak{(sl_n \times sl_m)}$-module isomorphic to
$L^{(n)}(k\omega_1+l\omega_{n-1})\boxtimes
L^{(m)}(l'\omega_1+k'\omega_{m-1})$. Also, the $U_q \mathfrak{(sl_n \times
sl_m)}$-action does not increase $\lambda(f)$, so
$L^{(n)}(k\omega_1+l\omega_{n-1})\boxtimes
L^{(m)}(l'\omega_1+k'\omega_{m-1}) \subset
\mathscr{E}_{N-1}(\Xi_{n,m})_{q,a}$ if $k+l' \leq a$. The fact that a
direct sum of all such modules exhaust $\mathscr{E}_{N-1}(\Xi_{n,m})_{q,a}$
can be obtained by calculating the dimensions. In fact, we have to verify
that
\begin{equation*}
\mathrm{dim} \mathscr{E}_{N-1}(\Xi_{n,m})_{q,a} \leq \sum_{k+l' \leq a} \mathrm{dim} \left(
L^{(n)}(k\omega_1+l\omega_{n-1})\boxtimes L^{(m)}(l'\omega_1+k'\omega_{m-1})\right).
\end{equation*}
Since one has the relation $\xi_1=0$ in $\mathscr{E}(\Xi_{n,m})_q$, the dimension of
$\mathscr{E}_{N-1}(\Xi_{n,m})_{q,a}$ satisfies the following inequality:
\begin{equation*}
\mathrm{dim} \mathscr{E}_{N-1}(\Xi_{n,m})_{q,a} \leq (C_{a+N-1}^{N-1})^2.
\end{equation*}
It is sufficient to verify the inequality
\begin{equation*}
(C_{a+N-1}^{N-1})^2 \leq \sum_{k+l' \leq a} \mathrm{dim} \left( L^{(n)}(k\omega_1+l\omega_{n-1})\boxtimes
L^{(m)}(l'\omega_1+k'\omega_{m-1})\right)
\end{equation*}
in the classical case. In the classical context this can be obtained via an
induction argument in $a$. \hfill $\square$

\bigskip

We are going to establish the necessary conditions for $\mathscr{E}_{s}(\Xi_{n,m})_q$ to be equivalent as
$U_q\mathfrak{sl}_N$-modules.

A special construction associates to every finite dimensional representation $V$ of $U_q\mathfrak{sl}_N$ a
central element $C_V$ of some extended algebra $U_q^\mathrm{ext}\mathfrak{sl}_N\supset U_q\mathfrak{sl}_N$
\cite{Kl-Sch}. It follows that the collection of constants $C_{L(\omega_p)}$, $p=1,\ldots,N$, constitute an
invariant of isomorphism for $\mathscr{E}_{2s}(\Xi_{n,m})_q$ as $U_q\mathfrak{sl}_N$-modules.

An essential property of the elements $C_V$ is that their action on the Verma module $M(\lambda)$ with
highest weight $\lambda$ is given by the constant \cite{Drin} (see also \cite[Proposition 3.1.22]{Vaksman}
for the special case $q \in (0,1)$)
$$C_V|_{M(\lambda)}=\sum_{\mu\in P}(\dim V_\mu)q^{-2(\mu,\lambda+\rho)},$$
where $P$ is the weight lattice of $U_q\mathfrak{sl}_N$, $V_\mu$ is the subspace of $\mu$-weight vectors in $V$, and $\rho$
is the half-sum of positive roots of $U_q\mathfrak{sl}_N$. Hence the same formula is applicable to any highest weight
$U_q\mathfrak{sl}_N$-module with highest weight $\lambda$.

A routine verification that involves \eqref{act_on_t_i}, \eqref{act_on_t_i*}, and \eqref{uq_acts_on_f(xin+1)}
shows that for $s=k \in \mathbb Z_+$  the vectors $t_1^kt_N^{*k} \in \mathscr{E}_{2s}(\Xi_{n,m})_q$ are also
$U_q\mathfrak{sl}_N$-singular (annihilated by $E_n$), and thus generate simple
$U_q\mathfrak{sl}_N$-submodules with highest weights $k(\omega_1+\omega_{N-1})$ for all $k\in\mathbb{Z}_+$.

A direct computation of those constants provides the result as follows. Let $e_p$ be the elementary symmetric degree $p$
polynomial in $N$ variables. Then
$$
C_{L(\omega_p)}|_{M(k(\omega_1+\omega_{N-1}))}=
e_p\left(q^{-2k-N+1},q^{-N+3},q^{-N+5},\ldots,q^{N-5},q^{N-3},q^{2k+N-1}
\right).
$$

 On the other hand, it is clearly visible from the definitions that the matrix elements of $U_q\mathfrak{sl}_N$-actions in
 $\mathscr{E}_{2k}(\Xi_{n,m})_q$ with respect to a suitable PBW-basis are Laurent polynomials of $q^{2k}$.

Thus an analytic continuation argument implies that the collection of constants
$$
e_p\left(q^{-2s-N+1},q^{-N+3},q^{-N+5},\ldots,q^{N-5},q^{N-3},q^{2s+N-1}
\right),\qquad p=1,\ldots,N,
$$
realizes as action of the central elements $C_{L(\omega_p)}$ on some non-zero simple submodules of the
$U_q\mathfrak{sl}_N$-modules $\mathscr{E}_{2s}(\Xi_{n,m})_q$.

Hence for isomorphic $\mathscr{E}_{2s}(\Xi_{n,m})_q$ and $\mathscr{E}_{2s'}(\Xi_{n,m})_q$, one should have the latter
collection of constants identical. This already implies that the collection of constants
$$q^{-2s-N+1},q^{-N+3},q^{-N+5},\ldots,q^{N-5},q^{N-3},q^{2s+N-1}$$
must be also identical, which means that, given such pair $s,s'$ then either $q^{-2s-N+1}=q^{-2s'-N+1}$ or
$q^{-2s-N+1}=q^{2s'+N-1}$. We obtain

\begin{proposition}
Given $s\in\mathbb{C}$, the set of those $s'$ for which $\mathscr{E}_{2s'}(\Xi_{n,m})_q$ is isomorphic to
$\mathscr{E}_{2s}(\Xi_{n,m})_q$ as $U_q\mathfrak{sl}_N$-modules, is a subset of
$$
\left\{\left. s+\frac{\pi in}{\ln q}\right|\:n\in\mathbb{Z}\right\}\cup\left\{\left. -(s+N-1)+\frac{\pi
in}{\ln q} \right|\:n\in\mathbb{Z}\right\}.
$$
\end{proposition}


\end{document}